\documentclass[10pt,reqno]{amsart}

\usepackage[a4paper,margin=27mm]{geometry}

\usepackage{xcolor}

\usepackage{bm}
\usepackage{amsmath,amssymb}
\usepackage{amsthm}
\usepackage{eurosym}
\usepackage{graphicx}
\usepackage{epsfig}
\usepackage{float}
\usepackage{hyperref}
\usepackage{dsfont}
\usepackage{booktabs}
\usepackage[nocompress, space]{cite}

\usepackage{enumitem}

\usepackage[mathlines]{lineno}
\usepackage{lineno}
\allowdisplaybreaks

\usepackage{ifthen}

\def\leq{\leqslant}
\def\geq{\geqslant}

\numberwithin{equation}{section}

\newtheoremstyle{thmlemcorr}{10pt}{10pt}{\itshape}{}{\bfseries}{.}{10pt}{{\thmname{#1}\thmnumber{
                        #2}\thmnote{ (#3)}}}
\newtheoremstyle{thmlemcorr*}{10pt}{10pt}{\itshape}{}{\bfseries}{.}\newline{{\thmname{#1}\thmnumber{
                        #2}\thmnote{ (#3)}}}
\newtheoremstyle{defi}{10pt}{10pt}{\itshape}{}{\bfseries}{.}{10pt}{{\thmname{#1}\thmnumber{
                        #2}\thmnote{ (#3)}}}
\newtheoremstyle{remexample}{10pt}{10pt}{}{}{\bfseries}{.}{10pt}{{\thmname{#1}\thmnumber{
                        #2}\thmnote{ (#3)}}}
\newtheoremstyle{ass}{10pt}{10pt}{}{}{\bfseries}{.}{10pt}{{\thmname{#1}\thmnumber{
                        A#2}\thmnote{ (#3)}}}

\theoremstyle{thmlemcorr}
\newtheorem{theorem}{Theorem}
\numberwithin{theorem}{section}

\theoremstyle{thmlemcorr*}
\newtheorem{theorem*}{Theorem}
\newtheorem{lemma*}[theorem]{Lemma}
\newtheorem{corollary*}[theorem]{Corollary}
\newtheorem{proposition*}[theorem]{Proposition}
\newtheorem{problem*}[theorem]{Problem}
\newtheorem{conjecture*}[theorem]{Conjecture}

\theoremstyle{defi}

\theoremstyle{remexample}
\newtheorem{remark}[theorem]{Remark}

\def\BibTeX{{\rm B\kern-.05em{\sc i\kern-.025em b}\kern-.08em
    T\kern-.1667em\lower.7ex\hbox{E}\kern-.125emX}}

\def\leq{\leqslant}
\def\geq{\geqslant}

\newtheoremstyle{thmlemcorr}{8pt}{8pt}{\itshape}{}{{\bfseries}}{.}{8pt}
  {{\thmname{#1}\thmnumber{ #2}\thmnote{ (#3)}}}
\newtheoremstyle{remexample}{8pt}{8pt}{}{}{{\bfseries}}{.}{8pt}
  {{\thmname{#1}\thmnumber{ #2}\thmnote{ (#3)}}}

\numberwithin{equation}{section}

\hypersetup{colorlinks=true,linkcolor=black,citecolor=black,urlcolor=black}
\AtBeginDocument{}

\definecolor{issueblue}{RGB}{0,55,175}
\definecolor{suggestionred}{RGB}{170,0,0}
\newcommand{\criticaltext}[1]{\textcolor{issueblue}{#1}}

\usepackage{placeins}

\begin{document}

\title{Beyond Periodicity: Fourier Methods for Optimal Control}

\author[Milena Ispiryan]{Milena Ispiryan}
\address{DiGREA Lab, CSIE, National Polytechnic University of Armenia}
\author[Vahagn Makaryan]{Vahagn Makaryan}
\address{DiGREA Lab, CSIE, Institute of Mathematics SNCO}
\author[Lusine Poghosyan]{Lusine Poghosyan}
\address{DiGREA Lab, CSIE, Institute of Mathematics SNCO, Armenia; Yerevan State University}
\author[Tigran Bakaryan]{Tigran Bakaryan}
\address{DiGREA Lab, CSIE, Institute of Mathematics SNCO, Armenia; Yerevan State University}
\email{tigran.bakaryan@instmath.sci.am}
\thanks{The research is supported by the Higher Education and Science Committee of the MESCS RA under Research Project No. 24IRF-1A001.}
\hypersetup{pdfauthor={Milena Ispiryan, Vahagn Makaryan, Lusine Poghosyan, Tigran Bakaryan},pdftitle={Beyond Periodicity: Fourier Methods for Optimal Control}}

\begin{abstract}
Fourier-based spectral methods provide compact global approximations, permit simple differentiation and integration, and achieve spectral convergence for smooth periodic functions. Their use in finite-horizon optimal control problems (OCPs) is nevertheless challenged by the generally nonperiodic state and control trajectories: endpoint mismatch induces artificial discontinuities and can destroy spectral convergence. Leveraging Eckhoff-type endpoint-correction techniques, we develop a Fourier-based spectral framework for OCPs with possibly nonlinear dynamics and constraints. The correction improves convergence without enlarging the set of decision variables, yielding an efficient finite-dimensional numerical method. Furthermore, we formulate selected dynamics and obstacle-avoidance
constraints through Sobolev residual bounds, which incorporate
information about the residuals over the entire time horizon
and improve constraint satisfaction. Numerical experiments on benchmark problems demonstrate competitive computational performance, lower residual errors, and improved solution accuracy.
\end{abstract}

\maketitle

\section{Introduction}
\label{sec:introduction}

 Optimal control problems arise in a wide range of applications,
including robotics, aerospace systems, process control,
and trajectory planning \cite{Betts2010, LaValle2006}. 
In many applications, the state and control trajectories
must satisfy dynamical equations and path constraints
over the entire time horizon. For problems with sufficiently
smooth trajectories, spectral methods provide an attractive
numerical approach, allowing accurate approximation using
relatively few global basis functions. By representing the
state and control trajectories in a finite-dimensional
basis, the original optimal control problem is approximated
by an optimization problem in the corresponding coefficients
\cite{Betts2010, Boyd2001}.
Among spectral methods, Fourier approximations play an
important role due to their compact representation and
favorable convergence properties for smooth periodic
functions. However, these properties are generally lost
when the endpoint compatibility conditions are not satisfied.
In this paper, leveraging Eckhoff's endpoint-correction
technique, we develop a Fourier-based approximation method
for optimal control problems with nonperiodic state and
control trajectories.

One of the most established approaches to numerical optimal
control are polynomial-based spectral methods, particularly
those using Chebyshev and Legendre polynomials.
These methods provide high-order approximation accuracy for
sufficiently smooth trajectories without requiring periodic
endpoint conditions. Moreover, their polynomial structure
allows accurate differentiation and integration using spectral
differentiation and quadrature techniques. The use of
Gauss--Lobatto nodes, which include both endpoints, also
facilitates the enforcement of initial and terminal conditions
\cite{Boyd2001,Garg2010}.Despite these advantages, Chebyshev and Legendre
 methods typically employ nonuniform
collocation grids whose nodes depend on the polynomial
degree. This can be inconvenient in applications requiring
fixed uniform sampling times, such as the sensing problem
considered in Section~\ref{subsec:fixed_rate_sensing}.
Moreover, approximating highly oscillatory trajectories
may require relatively high polynomial degrees (see Section~\ref{subsec:bicycle_high_frequency}), increasing
the number of optimization variables.
Polynomial spectral differentiation can also amplify
errors, with worst-case first-derivative amplification
growing quadratically with the polynomial degree,
compared with linear growth for Fourier differentiation
\cite{Boyd2001}.
Fourier approximations therefore provide an attractive
alternative for oscillatory trajectories, offering a
compact harmonic representation, simple differentiation,
and convenient evaluation on uniform grids. These advantages motivate our work on extending Fourier-based
optimal control methods to nonperiodic settings.

Bernstein polynomial approximations provide another
established approach to optimal control and trajectory
planning. Their convex-hull property and simple endpoint
evaluation facilitate the treatment of boundary and
trajectory constraints \cite{Farouki2012, Cichella2021, Cichella2022}.
Moreover, for sufficiently smooth functions, classical
Bernstein approximations converge uniformly together with
their derivatives up to any prescribed finite order
\cite{Cichella2022}.
However, the classical Bernstein approximation operator
has a relatively slow convergence rate, which may require
high polynomial degrees to accurately represent complex
and oscillatory trajectories. 
To combine the advantages of Bernstein polynomials and
Fourier approximations, mixed Bernstein--Fourier methods
have recently been proposed, using Bernstein polynomials
to represent nonperiodic behavior and Fourier modes to
capture periodic components \cite{Mudrik2025}.
In contrast, the Fourier--Eckhoff approach considered
in this paper accounts for nonperiodic behavior through
endpoint corrections without introducing additional
polynomial coefficients. The numerical experiments presented in this paper
further illustrate the advantages and limitations
of these approaches in terms of approximation accuracy,
constraint satisfaction, and computational efficiency.

Despite these advantages, the convergence properties of
Fourier approximations depend strongly on the periodic
compatibility of the approximated functions.
In finite-horizon optimal control problems, state and control
trajectories are generally nonperiodic, even when they are
smooth on the entire time interval.
Mismatches in endpoint values introduce artificial
discontinuities in their periodic extensions, leading to
the Gibbs phenomenon and preventing uniform convergence.
Moreover, mismatches in higher-order derivatives reduce
the decay rate of Fourier coefficients and consequently
slow down convergence.
These limitations motivate endpoint-correction techniques,
which account for the nonperiodic behavior and improve
the convergence properties of Fourier approximations.

To overcome the limitations of Fourier approximation
caused by endpoint incompatibility, several endpoint-correction
techniques have been developed for approximating nonperiodic
functions. These methods separate the nonperiodic behavior
from the part represented by a Fourier series.
The classical Krylov--Lanczos procedure constructs a
correction polynomial using known endpoint values and
derivatives, yielding a remainder with improved periodic
compatibility \cite{krylov1907approximate,lanczos1966discourse}.
Unlike the Krylov--Lanczos approach, Eckhoff's method
does not require prior knowledge of the endpoint values
or derivatives. Instead, it estimates the unknown endpoint
jumps from the asymptotic behavior of high-frequency
Fourier coefficients
\cite{Eckhoff1995, barkhudaryan2007asymptotic}.
We review these two approaches in
Sections~\ref{subsec:Krylov-Lanczos}
and~\ref{subsec:eckhoff_prelim}, respectively,
and use them as the foundation for the Fourier-based
optimal control formulation developed in this paper.

Leveraging these approximation techniques, we develop
a Fourier--Eckhoff method for finite-horizon optimal
control problems with possibly nonlinear dynamics and
constraints. The main idea is to parameterize the state
and control trajectories by Fourier coefficients and
express the unknown endpoint jumps directly in terms
of these coefficients. Unlike the classical Eckhoff
reconstruction, where the Fourier coefficients of a
given function are known, in our formulation they are
optimization variables. The endpoint corrections are
computed as linear functions of these variables and
therefore do not introduce additional independent
optimization variables.
The resulting approximation extends the applicability
of Fourier methods to nonperiodic settings while retaining
their main computational advantages, including compact
representation of oscillatory trajectories, analytical
differentiation, and convenient evaluation on uniform
time grids. These properties are particularly useful
for applications involving high-frequency oscillations
and fixed-rate measurements, as illustrated in the
numerical experiments in
Section~\ref{sec:numerical_results}.
The proposed formulation is developed in
Section~\ref{subsec:FE-OPC}.

To obtain a fully finite-dimensional optimal control problem,
the dynamics and path constraints are typically enforced
at a finite set of discretization nodes. However, satisfying
these constraints at the selected nodes does not guarantee
their satisfaction over the entire time interval.
In particular, constraint violations may occur between
the nodes, which can be problematic in safety-critical
applications. To address this limitation, we introduce
a Sobolev residual formulation for enforcing dynamics
and obstacle-avoidance constraints.
Unlike standard pointwise enforcement, the proposed
approach incorporates information about the residuals
over the entire time interval through integral constraints.
Using Sobolev embedding, we establish conditions under
which these constraints provide uniform residual bounds.
The numerical implementation uses discrete approximations
of the Sobolev norms, and constraint satisfaction is
evaluated on independent validation grids.
This extension, including the use of slack trajectories
for obstacle avoidance, is presented in
Section~\ref{subsec:sobolev_obstacle_avoidance}.
The numerical results demonstrate improved dynamics
residual control and positive sampled safety margins,
whereas the tested pointwise methods exhibit small
violations of the prescribed safety margins.
Moreover, increasing the approximation order of these
pointwise methods does not necessarily eliminate such
violations.

The performance of the proposed methods is investigated
through several numerical experiments.
We first consider a minimum-energy optimal control problem
with an explicit solution to evaluate approximation errors
and computational time
(Section~\ref{subsec:explicit_validation}).
We then examine high-frequency trajectory tracking with
a kinematic bicycle model
(Section~\ref{subsec:bicycle_high_frequency})
and a bearing-only localization problem with fixed-rate
measurements
(Section~\ref{subsec:fixed_rate_sensing}).
Finally, we investigate the effectiveness of Sobolev
residual constraints in a five-obstacle avoidance problem
(Section~\ref{subsec:sobolev_obstacle_avoidance}).
The proposed methods are compared with discrete,
Bernstein, Chebyshev, Legendre, Fourier--polynomial
\cite{NagurkaYen1990}, and mixed Bernstein--Fourier
\cite{Mudrik2025} approximations.
The comparisons consider approximation accuracy,
objective values, constraint satisfaction, number of
optimization variables, and computational time.
The results illustrate the advantages and limitations
of the different approaches, particularly for
oscillatory trajectories, fixed-rate sensing, and
safety-critical applications.

\section{Preliminaries}
\label{sec:preliminaries}

In this section, we recall classical results on Fourier approximation of functions defined on finite intervals, emphasizing how smoothness and endpoint compatibility affect convergence rates. We then present the Krylov--Lanczos endpoint-correction procedure and Eckhoff's asymptotic jump-recovery method, which account for endpoint mismatches and improve the approximation of nonperiodic functions. Finally, we introduce the class of finite-horizon nonlinear optimal control problems considered in this paper.

\subsection{Fourier-Based Function Approximation}
\label{subsec:functionapproximation}

We begin by recalling the classical truncated Fourier representation on a finite interval, which serves as the foundation for the endpoint-corrected constructions introduced in the subsequent subsections.

The truncated Fourier approximation of order $N$ of a function $f\in L^2[0,T]$ is defined by
\begin{equation}
\label{eq:prelim_fourier_series}
f_N(t) = \frac{a_0}{2} + \sum_{n=1}^{N} \Bigl( a_n\cos(\omega_n t) + b_n\sin(\omega_n t) \Bigr), \quad t\in[0,T],
\end{equation}
where the harmonic frequencies are explicitly given by
\begin{equation}
\label{eq:prelim_fourier_frequency}
\omega_n=\frac{2\pi n}{T}, \quad n=0,1,\dots,N.
\end{equation}
The associated Fourier coefficients are uniquely determined by orthogonal projection onto the finite-dimensional trigonometric subspace according to the classical Euler--Fourier formulas:
\begin{align}
\label{eq:prelim_an}
a_n &= \frac{2}{T}\int_0^T f(t)\cos(\omega_n t)\,dt, \quad n=0,1,\dots,N, \\
\label{eq:prelim_bn}
b_n &= \frac{2}{T}\int_0^T f(t)\sin(\omega_n t)\,dt, \quad n=1,2,\dots,N.
\end{align}

The convergence rate of the Fourier approximation is determined jointly by the smoothness of the target function on \([0,T]\) and by the extent to which its values and derivatives agree at the two endpoints. Specifically, let \(f\in C^{q+1}([0,T])\), \(q\geq 1\) and suppose that
\begin{equation}\label{eq-comp-cond}
f^{(j)}(0)=f^{(j)}(T), \quad j=0,1,\dots,q-1.
\end{equation}
Assume further that the $q$-th derivative does not
satisfy the endpoint compatibility condition; that is,
\begin{equation*}
f^{(q)}(0)\neq f^{(q)}(T).
\end{equation*}
Then, repeated integration by parts yields the following algebraic decay estimates for the Fourier coefficients (see \cite{barkhudaryan2007asymptotic})
\begin{equation}
\label{eq:alg_coeff_decay}
|a_n|,\ |b_n| = \mathcal{O}\left(\frac{1}{n^{q+1}}\right) \quad \text{as } n\to\infty.
\end{equation}
In turn, the global approximation error measured in the $L^2[0,T]$ norm satisfies the asymptotic rate (see \cite{barkhudaryan2007asymptotic})
\begin{equation}
\label{eq:l2_error_rate}
\|f-f_N\|_{L^2[0,T]} = \mathcal{O}\left(\frac{1}{N^{q+1/2}}\right),
\end{equation}
whereas the uniform error on the compact interval satisfies
\begin{equation}
\label{eq:linf_error_rate}
\|f-f_N\|_{\infty} = \max_{t\in[0,T]}|f(t)-f_N(t)| = \mathcal{O}\left(\frac{1}{N^q}\right),
\end{equation}
see \cite{zygmund2002trigonometric}. Note that when $q=0$, the $L^2$ convergence estimate, see \eqref{eq:l2_error_rate}
remains valid.
However, uniform convergence on $[0,T]$ fails due to
the endpoint discontinuity in the periodic extension
and the associated Gibbs phenomenon.
Nevertheless, pointwise convergence still holds
for every $t\in(0,T)$, see \cite{zygmund2002trigonometric}.

These estimates show that the convergence rate depends on the order at which endpoint compatibility, \eqref{eq-comp-cond}, fails. Thus, even when a function is smooth on the entire interval, mismatches in its endpoint values or derivatives can significantly reduce the convergence rate of its Fourier approximation. This limitation motivates the Krylov--Lanczos construction (see \cite{krylov1907approximate, lanczos1966discourse}), in which a correction term is chosen to reproduce the endpoint differences of the function and its derivatives. 

\subsection{Krylov-Lanczos Endpoint-Corrected Approximation}
\label{subsec:Krylov-Lanczos}
 When the endpoint compatibility conditions in \eqref{eq-comp-cond} are not satisfied, the Krylov--Lanczos procedure \cite{krylov1907approximate, lanczos1966discourse} improves the Fourier approximation by introducing a correction polynomial that accounts for the endpoint mismatches. The main assumption is that the values of the function and its derivatives at both endpoints, up to a prescribed order, are known. The correction polynomial is constructed using Bernoulli polynomials, and subtracting it from the original function yields a remainder whose values and derivatives satisfy the endpoint compatibility conditions up to the prescribed order, leading to an improved convergence rate of the Fourier approximation.

For a prescribed correction order \(q \geq 1\), denote the endpoint jumps
\begin{equation}
\label{eq:kl_jump_def}
A_k=
f^{(k)}(T)-f^{(k)}(0),
\qquad
k=0,1,\ldots,q-1,
\end{equation}
and assume that these quantities are known. The Krylov--Lanczos construction uses a standard family of correction functions \(\{B_k(t)\}_{k=0}^{q-1}\), given by suitably normalized Bernoulli polynomials, to construct a correction polynomial as a finite linear combination of these functions. The defining property of these functions is that each \(B_k\) carries exactly one unit endpoint jump in the corresponding derivative order:
\begin{equation}
\label{eq:kl_basis_jump_property}
B_k^{(j)}(T)-B_k^{(j)}(0)=
\delta_{kj},
\qquad
j,k=0,1,\ldots,q-1,
\end{equation}
where $\delta_{kj}$ denotes the Kronecker delta. Thus, each correction function $B_k$ produces a unit jump in its $k$-th derivative and zero jumps in all other derivatives included in the correction. 

The Bernoulli correction functions $B_k$ are defined by
\begin{equation}
\label{eq:kl_basis_B0}
B_0(t)=
\frac{t}{T}-\frac{1}{2},
\end{equation}
and, for $k=1,\ldots,q-1$, recursively by
\begin{equation}
\label{eq:kl_basis_recursive}
\frac{d}{dt}B_k(t)=
B_{k-1}(t),
\end{equation}
together with the zero-mean normalization
\begin{equation}
\label{eq:kl_basis_normalization}
\int_0^T B_k(t)\,dt=
0,
\qquad
k=0,1,\ldots,q-1.
\end{equation}
For example, the first two correction functions generated by this recursion and normalization are
\begin{equation*}
B_1(t)
=
\frac{t^2}{2T}
-
\frac{t}{2}
+
\frac{T}{12},\quad B_2(t)
=
\frac{t^3}{6T}
-
\frac{t^2}{4}
+
\frac{Tt}{12}.
\end{equation*}
The truncated Fourier approximation of order \(N\) of \(B_{k}(t)\) is given by
\begin{equation}
\label{eq:Bernoulli_real_Fourier}
B_{k,N}(t)=
\sum_{n=1}^{N}
\left(
\alpha_{k,n}\cos(\omega_n t)
+
\beta_{k,n}\sin(\omega_n t)
\right),
\end{equation}
where the cosine coefficients are
\begin{equation}\label{B_kcosine}
\alpha_{k,n}=\begin{cases}
0,
& k \text{ is even},\\[2mm]
\displaystyle
\frac{(-1)^{(k-1)/2}T^k}
{2^k\pi^{k+1}n^{k+1}},
& k \text{ is odd},
\end{cases}
\end{equation}
and the sine coefficients are
\begin{equation}\label{B_ksine}
\beta_{k,n}=\begin{cases}
\displaystyle
\frac{(-1)^{k/2+1}T^k}
{2^k\pi^{k+1}n^{k+1}},
& k \text{ is even},\\[2mm]
0,
& k \text{ is odd}.
\end{cases}
\end{equation}

Using these basis elements, we define the finite endpoint-correction polynomial of order $q$ as:
\begin{equation}
\label{eq:kl_poly_corr}
P_q(t)=\sum_{k=0}^{q-1} A_k B_k(t).
\end{equation}

The original function $f$ is decomposed into a correction polynomial and a
remainder 
\begin{equation}
\label{eq:kl_decomposition}
F(t)=f(t)-P_q(t) = f(t)-\sum_{k=0}^{q-1} A_k B_k(t),
\end{equation}
 with the same interior smoothness as $f$.
Note that by \eqref{eq:kl_jump_def} and \eqref{eq:kl_basis_jump_property}, the remainder $F(t)$ is guaranteed to meet the compatibility conditions, \eqref{eq-comp-cond}, up to the $(q-1)$-th derivative order:
\begin{equation}
\label{eq:kl_periodic_remainder}
F^{(k)}(T)=F^{(k)}(0), \qquad k=0,1,\dots,q-1.
\end{equation}

The improved periodic compatibility of $F$ permits the Fourier error estimates
above (see \eqref{eq:l2_error_rate} and  \eqref{eq:linf_error_rate}) to be applied at the corrected regularity order. Its Fourier projection is
\begin{equation}
\label{eq:kl_fourier_remainder}
F_N(t) = \frac{c_0}{2} +\sum_{n=1}^{N} \left( c_n\cos(\omega_n t)+d_n\sin(\omega_n t) \right),
\end{equation}
where the  coefficients are 
\begin{align}
\label{eq:kl_cn}
c_0 = a_0,\quad c_n =  a_{n}-\sum_{k=0}^{q-1}A_{k}\alpha_{k,n}, \quad
d_n = b_{n}-\sum_{k=0}^{q-1}A_{k}\beta_{k,n}, \quad n = 1, \dots, N,
\end{align}
where $\alpha_{k,n}$ and $\beta_{k,n}$ are given by
\eqref{B_kcosine} and \eqref{B_ksine}, respectively. 
Thus, the  Krylov--Lanczos approximation is then reconstructed by re-introducing the polynomial baseline
\begin{equation}
\label{eq:kl_final_approx}
I_{N,q}^{\mathrm{KL}}(t) = F_N(t)+\sum_{k=0}^{q-1}A_k B_k(t).
\end{equation}
 Because the corrected remainder satisfies the endpoint compatibility
conditions \eqref{eq-comp-cond}, the Krylov--Lanczos approximation
achieves the corresponding Fourier convergence rates, see for example, \cite{krylov1907approximate}.
More precisely, for $f\in C^{q+1}([0,T])$, with $q\geq 1$, we have
\begin{equation}
\label{eq:kl_l2_error}
\left\|
f
-
I_{N,q}^{\mathrm{KL}}
\right\|_{L^2[0,T]}=
\mathcal{O}\left(
\frac{1}{N^{q+\frac{1}{2}}}
\right),
\qquad
N\to\infty,
\end{equation}
and
\begin{equation}
\label{eq:kl_linf_error}
\max_{t\in[0,T]}
\left|
f(t)-I_{N,q}^{\mathrm{KL}}(t)
\right|=
\mathcal{O}\left(
\frac{1}{N^{q}}
\right),
\qquad
N\to\infty.
\end{equation}

\begin{remark}
Beyond the approximation of the trajectory functions themselves,
optimal control problems also require accurate approximations of
their derivatives. Here, we discuss the relationship between the
derivative of a truncated Fourier expansion and the truncated
Fourier expansion of the derivative, and examine how the
Krylov--Lanczos correction improves the approximation of derivatives.

Since the truncated Fourier expansion $S_N[f](t)=f_N(t)$
is a finite linear combination of trigonometric functions,
its derivative can be computed analytically by differentiating
term by term:
\begin{equation}
\label{eq:prelim_fourier_derivative}
\dot{f}_N(t)=
\sum_{n=1}^{N}
\omega_n
\Bigl(
-a_n(f)\sin(\omega_n t)
+b_n(f)\cos(\omega_n t)
\Bigr),
\qquad t\in[0,T].
\end{equation}

Suppose that $f\in C^1([0,T])$ satisfies the endpoint
compatibility condition
\begin{equation}
\label{eq:kl_F_periodic_zero}
f(T)=f(0).
\end{equation}
Then, differentiation commutes with the truncated Fourier
projection, that is,
\begin{equation*}
\frac{d}{dt}S_N[f](t)
=
\dot{f}_N(t)
=
S_N[\dot{f}](t).
\end{equation*}
Indeed, integration by parts shows that the Fourier
coefficients of $\dot{f}$ are $\omega_n b_n(f)$ and
$-\omega_n a_n(f)$, while $a_0(\dot{f})=0$.

However, when $f(T)\neq f(0)$, this property no longer
holds. In this case, integration by parts introduces
additional boundary terms, and differentiating the
truncated Fourier expansion of $f$ does not yield the
truncated Fourier expansion of $\dot{f}$. Hence, special care is needed when approximating derivatives.

In the case of Krylov--Lanczos approximation for the remainder $F=f-P_q$ (satisfies
the endpoint compatibility conditions up to order $q-1$).  Consequently, for $q\geq1$, we have
\begin{equation*}
\frac{d}{dt}S_N[F](t)
=
S_N[\dot{F}](t).
\end{equation*}
Differentiating the Krylov--Lanczos approximation
\eqref{eq:kl_final_approx}, we obtain
\begin{equation}
\label{eq:kl_final_derivative}
\dot{I}_{N,q}^{\mathrm{KL}}(t)=
\dot{F}_N(t)
+
\frac{A_0}{T}
+
\sum_{k=1}^{q-1}A_k B_{k-1}(t),
\end{equation}
where the last sum is understood to be empty when $q=1$.

We now state the corresponding convergence estimates for the derivatives. Let $f\in C^{q+1}([0,T])$,
with $q\geq2$. Since
\begin{equation}
\label{eq:kl_derivative_error_identity}
\dot{f}(t)
-
\dot{I}_{N,q}^{\mathrm{KL}}(t)
=
\dot{F}(t)-\dot{F}_N(t),
\end{equation}
and the corrected remainder satisfies
\begin{equation}
\label{eq:kl_derivative_periodicity}
\dot{F}^{(j)}(T)
=
\dot{F}^{(j)}(0),
\qquad
j=0,1,\ldots,q-2,
\end{equation}
the standard Fourier approximation estimates applied
to $\dot{F}\in C^q([0,T])$ yield
\begin{equation*}
\left\|
\dot{f}
-
\dot{I}_{N,q}^{\mathrm{KL}}
\right\|_{L^2[0,T]}
=
\mathcal{O}\left(
\frac{1}{N^{q-\frac12}}
\right),\qquad
\max_{t\in[0,T]}
\left|
\dot{f}(t)
-
\dot{I}_{N,q}^{\mathrm{KL}}(t)
\right|
=
\mathcal{O}\left(
\frac{1}{N^{q-1}}
\right),
\qquad N\to\infty.
\end{equation*}

Thus, the Krylov--Lanczos correction not only improves
the approximation of the original function but also
provides improved convergence rates for its derivatives.
\end{remark}

\subsection{Asymptotic Jump Estimation via Eckhoff's Method}
\label{subsec:eckhoff_prelim}

The Krylov--Lanczos correction introduced in the previous
subsection requires the endpoint jumps
$A_k=f^{(k)}(T)-f^{(k)}(0)$,
$k=0,1,\ldots,q-1$, to be known in advance.
In many applications, however, these quantities are unknown.
Eckhoff's method estimates them from the high-frequency
behavior of the Fourier coefficients of $f$;
see \cite{barkhudaryan2007asymptotic}.

For $f\in C^{q+1}([0,T])$, repeated integration by parts separates the endpoint
jumps from a remainder with coefficients of order $n^{-(q+1)}$, giving
the following  asymptotic expansions
\begin{align}
\label{eq:eckhoff_asymptotic_an}
a_n &= \sum_{k=0}^{q-1} A_k \alpha_{k,n}+\mathcal{O}\left(
\frac{1}{n^{q+1}}
\right), \qquad n \gg 1,\\
\label{eq:eckhoff_asymptotic_bn}
b_n &= \sum_{k=0}^{q-1} A_k \beta_{k,n}+\mathcal{O}\left(
\frac{1}{n^{q+1}}
\right), \qquad n \gg 1,
\end{align}
where \(\alpha_{k,n}\) and \(\beta_{k,n}\) are defined by \eqref{B_kcosine} and \eqref{B_ksine}, respectively. These relations show that the leading large-\(n\) behavior of the Fourier coefficients is determined by the endpoint jumps \(A_0,\ldots,A_{q-1}\). This observation provides the basis for estimating the jumps from the large-index Fourier coefficients.

To estimate the unknown endpoint jumps, define
\begin{equation}
\label{eq:eckhoff_mu_definition}
\mu_{k,n}
=
\alpha_{k,n}
+
\beta_{k,n}.
\end{equation}
Since \(\alpha_{k,n}=0\) for even \(k\) and \(\beta_{k,n}=0\) for odd \(k\), the coefficient \(\mu_{k,n}\) is nonzero for each \(k=0,1,\ldots,q-1\). Adding the asymptotic relations \eqref{eq:eckhoff_asymptotic_an} and \eqref{eq:eckhoff_asymptotic_bn} gives
\begin{equation}
\label{eq:eckhoff_asymptotic_combined}
a_n+b_n=
\sum_{k=0}^{q-1}
A_k\mu_{k,n}
+
\mathcal{O}\left(
\frac{1}{n^{q+1}}
\right),
\qquad
n\to\infty.
\end{equation}

We choose \(q\) distinct large indices
\[
n_0,n_1,\ldots,n_{q-1}
\]
from the available range of Fourier coefficients. Define

\begin{equation}\label{eq:eckhoff_tail_vector}
\mathbf{\Lambda}^{(N)}
=
\begin{bmatrix}
a_{n_0} + b_{n_0} \\
\vdots \\
a_{n_{q-1}} + b_{n_{q-1}}
\end{bmatrix},
\end{equation}

and
\begin{equation}
\label{eq:eckhoff_recovery_matrix}
\mathcal{B}^{(N)}
=
\begin{bmatrix}
\mu_{0,n_0} & \mu_{1,n_0} & \cdots & \mu_{q-1,n_0}\\
\mu_{0,n_1} & \mu_{1,n_1} & \cdots & \mu_{q-1,n_1}\\
\vdots & \vdots & \ddots & \vdots\\
\mu_{0,n_{q-1}} & \mu_{1,n_{q-1}} & \cdots & \mu_{q-1,n_{q-1}}
\end{bmatrix}.
\end{equation}
The vector of endpoint jumps is then estimated by solving the linear system
\begin{equation}
\label{eq:eckhoff_recovery_system}
\mathcal{B}^{(N)}
\widehat{\mathbf A}
=
\mathbf{\Lambda}^{(N)}.
\end{equation}

The matrix \(\mathcal{B}^{(N)}\) has a Vandermonde-type structure in the variables \(1/n_0,\ldots,1/n_{q-1}\). Hence, for distinct nonzero indices \(n_0,\ldots,n_{q-1}\), it is nonsingular in exact arithmetic. In practice, however, the system may become ill-conditioned when \(q\) is large or when the values \(1/n_0,\ldots,1/n_{q-1}\) are close to one another. When \(\mathcal{B}^{(N)}\) is nonsingular, the estimate can be written formally as
\begin{equation}
\label{eq:eckhoff_Ahat_explicit_square}
\widehat{\mathbf A}
=
\left(\mathcal{B}^{(N)}\right)^{-1}
\mathbf{\Lambda}^{(N)}.
\end{equation}

Using the estimated jump vector \(\widehat{\mathbf A}\), define the estimated endpoint-correction polynomial by
\begin{equation}
\label{eq:eckhoff_recovered_correction_real}
\widehat P_q(t)
=
\sum_{k=0}^{q-1}
\widehat A_k B_k(t)
=
\mathbf B(t)^\top \widehat{\mathbf A},
\end{equation}
where
\begin{equation}
\label{eq:eckhoff_B_vector}
\mathbf B(t)
=
\begin{bmatrix}
B_0(t) & B_1(t) & \cdots & B_{q-1}(t)
\end{bmatrix}^{\top}.
\end{equation}
Equivalently,

\begin{equation}
\label{eq:eckhoff_recovered_correction_matrix}
\widehat P_q(t)
=
\mathbf B(t)^\top
\left(\mathcal B^{(N)}\right)^{-1}
\mathbf{\Lambda}^{(N)}.
\end{equation}

We then define the corrected remainder by
\begin{equation}
\label{eq:eckhoff_corrected_remainder_definition}
\widehat F(t)
=
f(t)-\widehat P_q(t).
\end{equation}
The Fourier coefficients of this corrected remainder are computed from the original coefficients by subtracting the contribution of the estimated endpoint correction:
\begin{equation}\label{eq:eckhoff_corrected_cn}
\widehat c_n=
a_n-
\sum_{k=0}^{q-1}
\widehat A_k\alpha_{k,n},
\qquad
n=1,\ldots,N,
\end{equation}
\begin{equation}\label{eq:eckhoff_corrected_dn}
\widehat d_n=
b_n-
\sum_{k=0}^{q-1}
\widehat A_k\beta_{k,n},
\qquad
n=1,\ldots,N.
\end{equation}
Since the correction functions have zero mean, the constant coefficient is unchanged. Hence,
\begin{equation}
\label{eq:eckhoff_corrected_F_N}
\widehat F_N(t)
=
\frac{a_0}{2}
+
\sum_{n=1}^{N}
\left(
\widehat c_n\cos(\omega_n t)
+
\widehat d_n\sin(\omega_n t)
\right).
\end{equation}
The Eckhoff-corrected reconstruction is therefore
\begin{equation}
\label{eq:eckhoff_final_reconstruction_real}
\widetilde f_{N,q}(t)
=
\widehat F_N(t)
+
\widehat P_q(t)
=
\widehat F_N(t)
+
\mathbf B(t)^\top \widehat{\mathbf A}.
\end{equation}

The same Fourier coefficients \(a_n\) and \(b_n\) are used both to estimate the endpoint jumps and to construct the corrected residual coefficients. Thus, once the original Fourier coefficients have been computed, the correction step does not require any additional quadrature.

In this setting, Eckhoff's method can be viewed as a Krylov--Lanczos correction in which the endpoint jumps are estimated from the large-index Fourier coefficients rather than supplied in advance.

\subsection{Optimal Control Problem Formulation}
\label{subsec:ocp_formulation}

We now specify the class of finite-horizon optimal control problems to which the preceding approximation tools will be applied. The purpose of the Fourier, Krylov--Lanczos, and Eckhoff constructions is to provide finite-dimensional representations of trajectories that are defined on a finite time interval and need not satisfy periodic endpoint conditions. This is important in optimal control, where state and control trajectories are generally non-periodic and may therefore be poorly approximated by a direct truncated Fourier series near the endpoints.

Let \(T>0\) be fixed. We consider state trajectories
\[
x:[0,T]\to\mathbb{R}^{n_x}
\]
and control trajectories
\[
u:[0,T]\to\mathbb{R}^{n_u}
\]
satisfying a nonlinear dynamical system, endpoint constraints, and pointwise path constraints. The general problem is
\begin{equation}\label{eq:ocp_general}
\begin{aligned}
\min_{x,u}\quad &E(x(0),x(T))+
                   \int_0^T G(t,x(t),u(t))\,dt,\\
\text{subject to}\quad&\dot x(t)=f(t,x(t),u(t))\quad\text{a.e.},\\
&e(x(0),x(T))=0,\\
&h(t,x(t),u(t))\leq0.
\end{aligned}
\end{equation}
 The numerical dynamics in this paper are autonomous,
so the explicit argument $t$ can be omitted in those application formulas.
Tracking costs still depend on time through the prescribed reference.

Here \(x(t)\in\mathbb{R}^{n_x}\) denotes the state vector and \(u(t)\in\mathbb{R}^{n_u}\) denotes the control vector. The mapping
\[
E:\mathbb{R}^{n_x}\times\mathbb{R}^{n_x}\to\mathbb{R}
\]
defines the endpoint cost, while
$G:[0,T]\times\mathbb{R}^{n_x}\times\mathbb{R}^{n_u}\to\mathbb{R}$
defines the running cost. The dynamics are governed by
$f:[0,T]\times\mathbb{R}^{n_x}\times\mathbb{R}^{n_u}\to\mathbb{R}^{n_x}$.
The mappings
$e:\mathbb{R}^{n_x}\times\mathbb{R}^{n_x}\to\mathbb{R}^{n_e}$
and $h:[0,T]\times\mathbb{R}^{n_x}\times\mathbb{R}^{n_u}\to\mathbb{R}^{n_h}$
define the endpoint equality constraints and the pointwise path inequality constraints, respectively.

This formulation includes many standard finite-horizon trajectory optimization problems. For example, a quadratic tracking objective is obtained with
\begin{align*}
E(x(0),x(T))&=(x(T)-x_d)^\top S(x(T)-x_d),\\
G(t,x,u)&=(x-x_r(t))^\top Q(x-x_r(t))+u^\top Ru,
\end{align*}
where $Q,S\succeq0$ and $R\succ0$ are symmetric matrices of appropriate
dimensions. 
The fixed-rate sensing example (see Section \ref{subsec:fixed_rate_sensing}) later uses a sampled information functional
instead of this Bolza objective. Its dependence on finitely many trajectory
evaluations is handled by the same coefficient parameterization; it is an
explicit extension of the objective class in \eqref{eq:ocp_general}.

The problem \eqref{eq:ocp_general} is infinite-dimensional because the unknowns are the functions \(x(\cdot)\) and \(u(\cdot)\). The following sections use the endpoint-corrected Fourier representations developed above to approximate these trajectories and obtain a finite-dimensional optimization problem.

\section{Fourier-Based Approximations for Optimal Control Problems}
\label{sec:ocp_fourier}

The main approach considered in this paper for solving optimal
control problems is based on spectral methods. The idea is to
approximate the original infinite-dimensional problem by a
finite-dimensional optimization problem, where the decision
variables are the coefficients of the state and control
trajectories represented in a chosen basis of functions.
In this section, we present Fourier-based approximation methods,
with particular emphasis on the Krylov--Lanczos and
Eckhoff-corrected Fourier approximations introduced in
Section~\ref{sec:preliminaries}.

\subsection{Approximate State and Control Parameterization}
\label{subsec:general_parameterization}

Following the spectral approximation approach, let \(x_N:[0,T]\to\mathbb{R}^{n_x}\) and
\(u_N:[0,T]\to\mathbb{R}^{n_u}\) denote finite-dimensional approximations of
the state and control trajectories. We write
\begin{equation}
\label{eq:general_xN_uN}
x(t)\approx x_N(t),
\qquad
u(t)\approx u_N(t),
\end{equation}
where each component of \(x_N\) and \(u_N\) is represented by a Fourier-based
or mixed approximation chosen according to the qualitative structure of the
problem.
More precisely, for each state component \(x_i\), \(i=1,\dots,n_x\), and each
control component \(u_j\), \(j=1,\dots,n_u\), we introduce approximants of the
form
\begin{equation}
\label{eq:general_component_param}
x_{i,N}(t)=\Phi_i(t;\theta_{x_i}),
\qquad
u_{j,N}(t)=\Psi_j(t;\theta_{u_j}),
\end{equation}
where \(\theta_{x_i}\) and \(\theta_{u_j}\) are finite collections of
approximation coefficients. The complete decision vector is then
\begin{equation}
\label{eq:general_theta}
\theta=
\bigl(
\theta_{x_1},\dots,\theta_{x_{n_x}},
\theta_{u_1},\dots,\theta_{u_{n_u}}
\bigr).
\end{equation}

Substituting these approximants into the cost functional and constraints yields
the approximate problem
\begin{equation}
\label{eq:approx_ocp_cost}
\min_{\theta}\;
J_N(\theta)
=
E\bigl(x_N(0),x_N(T)\bigr)
+
\int_0^T G\bigl(t,x_N(t),u_N(t)\bigr)\,dt
\end{equation}
subject to
\begin{equation}
\label{eq:approx_ocp_dyn}
\left\|\dot x_N(t)-f\bigl(t,x_N(t),u_N(t)\bigr)\right\|_\infty\leq \delta,
\qquad t\in[0,T],
\end{equation}
\begin{equation}
\label{eq:approx_ocp_boundary}
e\bigl(x_N(0),x_N(T)\bigr)=0,
\end{equation}
and
\begin{equation}
\label{eq:approx_ocp_path}
h\bigl(t,x_N(t),u_N(t)\bigr)\le 0,
\qquad t\in[0,T].
\end{equation}
 In \eqref{eq:approx_ocp_dyn} choosing small enough $\delta$, we ensure 
$\dot x_N(t)-f\bigl(t,x_N(t),u_N(t)\bigr)\approx 0,
\quad t\in[0,T]$. The trajectory parameterization first gives a finite-dimensional but
semi-infinite problem: there are finitely many coefficients, but the path
conditions are imposed for every $t\in[0,T]$.

In numerical implementation, the dynamics and other path inequalities are enforced at specified nodes
(which may differ by constraint or method) making the problem finite dimensional. Particularly, the dynamic feasibility condition
\eqref{eq:approx_ocp_dyn} is enforced on a prescribed grid $0=t_0<t_1<\cdots<t_M=T$; that is,
\begin{equation}
\label{eq:approx_ocp_dyn_grid}
\left\|
\dot x_N(t_\ell)
-
f\bigl(t_\ell,x_N(t_\ell),u_N(t_\ell)\bigr)
\right\|_{\infty}
\le \delta,
\qquad \ell=0,1,\dots,M.
\end{equation}
where \(\delta\ge 0\) is a chosen tolerance.

The approximation framework is intentionally kept general at this stage. In the
subsections below, we specify three practically important cases: pure Fourier
approximation for periodic problems, Krylov--Lanczos correction for
non-periodic problems with known endpoint jump data, and Eckhoff correction for
non-periodic problems in which the jump data must be recovered from the
spectral coefficients.

\subsection*{Standard Fourier Approximation}

If the endpoint conditions admit periodic state and control approximants,
and the trajectories are sufficiently smooth under periodic continuation,
it is natural to use truncated Fourier series directly.
A standard Fourier expansion always has equal values at $0$ and $T$; being
``close to periodic'' cannot make incompatible hard endpoint conditions
feasible.

 In that case, one sets
\begin{align}
\label{eq:pure_fourier_x}
x_N(t)
&=
\frac{a_0^{x}}{2}
+
\sum_{n=1}^{N}
\Bigl(
a_n^{x}\cos(\omega_n t)+b_n^{x}\sin(\omega_n t)
\Bigr),\\
\label{eq:pure_fourier_u}
u_N(t)
&=
\frac{a_0^{u}}{2}
+
\sum_{n=1}^{N}
\Bigl(
a_n^{u}\cos(\omega_n t)+b_n^{u}\sin(\omega_n t)
\Bigr),
\end{align}
where \(\omega_n=2\pi n/T\), and the vectors
\(a_n^x,b_n^x\in\mathbb{R}^{n_x}\),
\(a_n^u,b_n^u\in\mathbb{R}^{n_u}\) are optimization variables.

This representation is compact and spectrally efficient when
the trajectories are smooth and satisfy periodic endpoint
compatibility conditions. However, in many finite-horizon
optimal control problems, these conditions are not satisfied,
reducing the effectiveness of standard Fourier approximations.
In particular, a mismatch in endpoint values prevents
uniform convergence on the entire interval and leads to
the Gibbs phenomenon near the endpoints.

\subsection{Fourier Approximation With Polynomial Endpoint Correction}
\label{subsec:poly_correction_prelim}

A classical way to reduce endpoint mismatch is to split the function into a
low-order polynomial part and a Fourier part. One writes
\begin{equation}
\label{eq:prelim_poly_fourier_split}
f(t)=P_m(t)+\lambda(t),
\end{equation}
where \(P_m\) is a polynomial of degree \(m\), chosen so as to satisfy
prescribed endpoint conditions, and \(\lambda\) is approximated by a truncated
Fourier series. 
This yields the approximation
\begin{equation}
\label{eq:prelim_poly_fourier_approx}
f_N^{\mathrm{poly}}(t)=P_m(t)+\lambda_N(t),
\end{equation}
The
remaining part is then represented by a Fourier series,
\begin{equation}
\label{eq:prelim_lambda_series}
\lambda_N(t)=
\sum_{n=1}^{N}
\Bigl(
\alpha_n\cos(\omega_n t)+\beta_n\sin(\omega_n t)
\Bigr). 
\end{equation}
To solve OPC problem, for the state and control are used this approximation.

In \cite{NagurkaYen1990}, a fifth-order
polynomial is used,
\begin{equation}
\label{eq:prelim_poly5}
P_5(t)=\sum_{j=0}^{5} p_j t^j.
\end{equation}
The six polynomial coefficients can be determined from prescribed values of
the total trajectory and its first two derivatives at $0$ and $T$, after
subtracting the corresponding Fourier contributions \cite{NagurkaYen1990}.  Any endpoint derivative not prescribed
by the OCP remains an optimization variable in the implemented baseline.
Thus, prescribing $P_5$ alone to have the desired endpoint values would
generally be insufficient: the Fourier part need not vanish there.

The polynomial correction approach is effective when endpoint values and
derivatives are available explicitly. However, the choice of polynomial degree
is fixed a priori and is not directly tied to the asymptotic jump structure of
the Fourier coefficients and it increase number of variables.

\subsection{Krylov--Lanczos  Approximation}\label{subsec:KL-OPC}

In some practical optimal control problems, the trajectories are non-periodic
but the endpoint data are sufficiently structured that the jump information can
be determined from the boundary conditions and smoothness requirements. Typical
examples include problems in which one prescribes
\[
x(0),\quad x(T),\quad \dot x(0),\quad \dot x(T),
\quad \ddot x(0),\quad \ddot x(T),
\]
or analogous endpoint conditions for the control.

In such cases, one may apply the Krylov--Lanczos correction component-wise. For
the available   endpoint jumps as before set  (coordinate-wise)
\begin{equation}
\label{eq:ocp_KL_jump}
A_k = x^{(k)}(T)-x^{(k)}(0),
\qquad k=0,1,\dots,q-1.
\end{equation} 
Then, the corrected approximation is 
\begin{equation}
\label{eq:ocp_KL_xN_scalar}
x_N(t)=F_N(t)+P_q(t),
\end{equation}
where
\begin{equation}
\label{eq:ocp_KL_Pq}
P_q(t)=\sum_{k=0}^{q-1}A_k B_k(t),
\end{equation}
and \(F_N\) is a truncated Fourier approximation vector of the corresponding periodic
remainder. Analogously, we construct an approximation of the control
trajectory using the $r$ available endpoint conditions.
In vector form, the state and control approximations are
written as
\begin{align}
\label{eq:ocp_KL_xN}
x_N(t)
&=
F_N^x(t)+P_q^x(t),\\
\label{eq:ocp_KL_uN}
u_N(t)
&=
F_N^u(t)+P_r^u(t).
\end{align}

The correction order is dictated by the amount of available endpoint
information. For example, if only the value mismatch is known, then a
first-order correction may be used; if both value and first-derivative jumps
are known, then a second-order correction becomes natural. More generally, the
Krylov--Lanczos correction is appropriate when the problem admits sufficiently
smooth feasible trajectories and the relevant endpoint jump quantities can be
computed or imposed consistently from the optimal control formulation.

\subsection{ Fourier-Eckhoff Approximation}\label{subsec:FE-OPC}
In this subsection, we present the main approach proposed in
this paper, the Fourier--Eckhoff approximation for optimal
control problems, which is applicable to both periodic and
nonperiodic settings. Unlike the Fourier approximation with
polynomial correction (see Section~\ref{subsec:poly_correction_prelim}), the proposed method does not introduce additional optimization
variables. Moreover, in contrast to the Krylov--Lanczos
approximation, it does not require prior knowledge of the
endpoint values or derivatives of the optimal trajectories.

In the standard Eckhoff setting, the function Fourier coefficients, $\mathbf{f}^{(N)}$, are known; that is, \eqref{eq:eckhoff_tail_vector}  is given. We choice correction order $q$ and $q$ distinct large indices
$n_0,n_1,\ldots,n_{q-1}
$. Then, leveraging \eqref{eq:eckhoff_asymptotic_combined}, we approximate the  endpoint jumps by solving following  linear system
\begin{equation}\label{eq:main}
\mathcal{B}^{(N)}
\widehat{\mathbf A}
=
\mathbf{f}^{(N)},
\end{equation}
where $\mathcal{B}^{(N)}$ is known (non-singular) and given  by \eqref{eq:eckhoff_tail_vector}. In the optimal control case, the Fourier  coefficients  of the  solution, $\mathbf{f}^{(N)}$ are unknown. Nevertheless, the estimate in \eqref{eq:eckhoff_asymptotic_combined} holds true for the solution, hence, the  we can use approximation in \eqref{eq:main} to express the endpoint jumps in terms of the optimization variables $\mathbf{f}^{(N)}$; that is,
\begin{equation}
\label{eq:ocp_eckhoff_Ahat}
\widehat{\mathbf A}
=
\left(\mathbf{\mathcal B}^{(N)}\right)^{-1}\mathbf f^{(N)}.
\end{equation}
Using this, we write the Eckhoff  correction term in terms of Fourier coefficients of the optimal solution
\begin{equation}
\label{eq:ocp_eckhoff_Pq}
\widehat P_q(t)
=
\sum_{k=0}^{q-1}\widehat A_k B_k(t) = \mathbf B(t)^\top \left(\mathbf{\mathcal B}^{(N)}\right)^{-1}\mathbf f^{(N)}.
\end{equation}

Therefore, the Eckhoff-corrected approximation becomes
\begin{equation}
\label{eq:ocp_eckhoff_xN_scalar}
x_N(t)
=
\widehat F_N(t)+\widehat P_q(t).
\end{equation}
 Repeating the above steps for the state and control, we get following approximations
\begin{align}
\label{eq:ocp_eckhoff_xN}
x_N(t)
&=
\widehat F_N^x(t)+\widehat P_q^x(t),\\
\label{eq:ocp_eckhoff_uN}
u_N(t)
&=
\widehat F_N^u(t)+\widehat P_r^u(t).
\end{align}

In the OCP, the Fourier coefficients are optimization variables rather than
known coefficients of an already available optimal trajectory.
The recovered jumps are determined from those variables by the fixed linear
map \criticaltext{$\mathcal{B}^{N}$}, and the corrected Fourier coefficients must be formed by
\eqref{eq:eckhoff_corrected_cn}--\eqref{eq:eckhoff_corrected_dn}.
Thus the correction is part of the trajectory parameterization, not an
additional set of free jump variables. For common order $N$, the state and
control representation has $(n_x+n_u)(2N+1)$ variables.

\section{Numerical Results and Applications}
\label{sec:numerical_results}

In this section, we present numerical experiments to evaluate
the performance of the proposed Fourier--Eckhoff approximation
for solving optimal control problems. We first validate the
method on a problem with a known optimal solution, followed
by high-frequency trajectory tracking and a fixed-rate
sensing application. The results are compared with alternative
numerical methods in terms of solution accuracy, constraint
satisfaction, and computational efficiency.

For an approximate state--control pair $(x_N,u_N)$,
we define the dynamics residual
\begin{equation}
r_N(t)=\dot{x}_N(t)-f(t, x_N(t),u_N(t)).
\label{eq:num_residual}
\end{equation}
Its reported maximum is evaluated on an independent
validation grid:
\begin{equation}
R_\infty
=
\max_{t\in\mathcal G_{\mathrm{val}}}
\|r_N(t)\|_\infty.
\label{eq:num_residual_max}
\end{equation}
Thus, $R_\infty$ is a sampled validation measure.
Successful solver termination alone does not establish
continuous-time feasibility.

\subsection{Validation against an Explicit Solution}
\label{subsec:explicit_validation}

We first consider an optimal control problem with a known
explicit solution to evaluate the accuracy of the proposed
Fourier--Eckhoff approximation. This allows us to directly
compare the numerical results with the exact optimal solution.
Specifically, we consider the following minimum-energy problem:
\begin{equation}
\begin{aligned}
\min_{x,u}\quad&
J(x,u)=\int_0^T u(t)^2\,dt,\\
\text{subject to}\quad&
\dot{x}(t)=x(t)+u(t),\\
&x(0)=2,\qquad x(T)=10,
\end{aligned}
\qquad T=2\pi.
\label{eq:num_explicit_problem}
\end{equation}
No additional state or control bounds are imposed.
\begin{figure}[htbp]
\centering
\includegraphics[width=\linewidth]{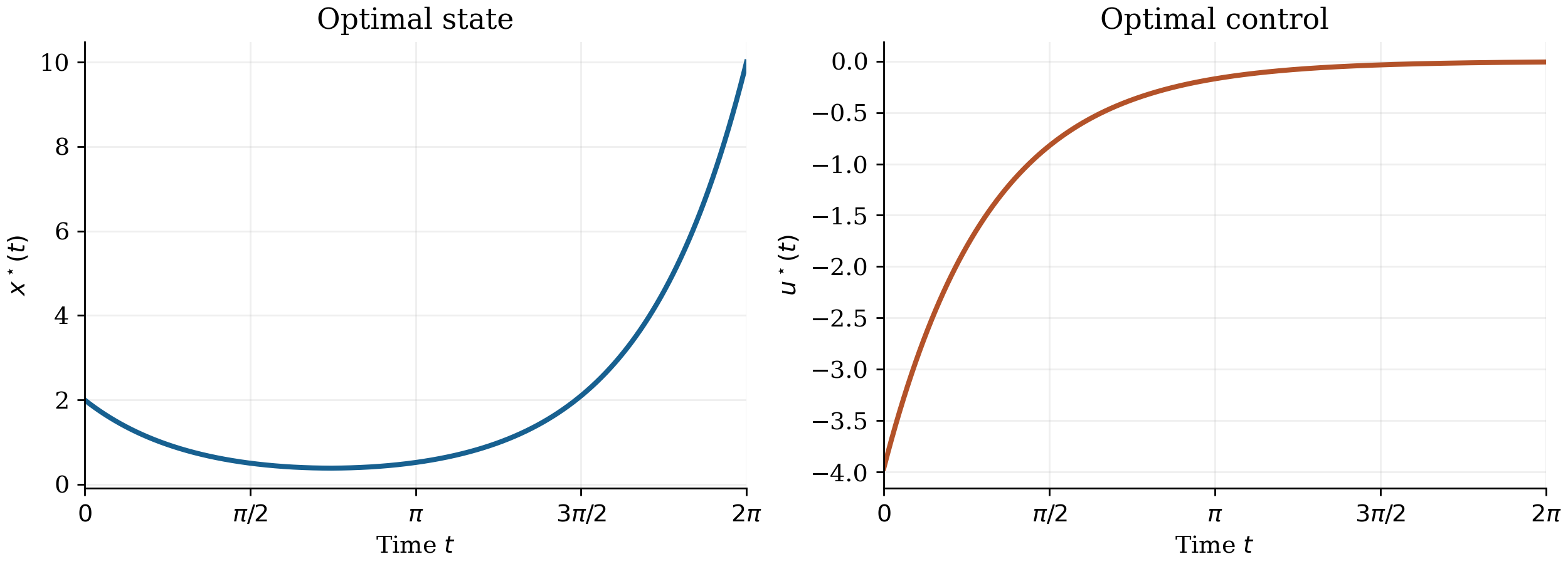}
\caption{Exact optimal state and control.}
\label{fig:explicit_energy_exact}
\end{figure}

The unique optimal state and control are
\begin{equation}
x^\star(t)=Ce^{-t}+De^{t-T},
\qquad
u^\star(t)=-2Ce^{-t},
\label{eq:num_explicit_solution}
\end{equation}
where
\begin{equation}
C=\frac{2-10e^{-T}}{1-e^{-2T}},
\qquad
D=\frac{10-2e^{-T}}{1-e^{-2T}}.
\end{equation}
Consequently,
\begin{equation}
J^\star
=
2C^2(1-e^{-2T})
\approx 7.8513294302.
\end{equation}

To see why this solution is optimal, write any other
feasible control as $u=u^\star+w$.
The terminal condition implies
$\int_0^T e^{-t}w(t)\,dt=0$, and therefore
\begin{equation}
J(u)-J(u^\star)=\int_0^T w(t)^2\,dt\geq0.
\end{equation}
Equality holds only when $w=0$ almost everywhere. Hence, $u^*$ is a minimizer.

The numerical formulations use $K=1001$ uniform
nodes and dynamics tolerance $\delta=10^{-6}$.
The polynomial methods use $d=32$.
FE and Krylov--Lanczos (KL) use $N=32$ and $q=2$;
KL is supplied with the exact endpoint jumps.
The Fourier--polynomial method uses $15$ Fourier modes,
a quintic state correction, and a linear control
correction. The mixed method uses Bernstein degree
$5$ and $13$ Fourier modes per signal.

On an independent grid of $16\,385$ points, we measure
\begin{align}
E_{x,\infty}
=
\max_{t\in\mathcal G_{\mathrm{val}}}
|x_N(t)-x^\star(t)|,\qquad
E_{u,\infty}
=
\max_{t\in\mathcal G_{\mathrm{val}}}
|u_N(t)-u^\star(t)|.
\end{align}

\begin{table}[htbp]
\centering
\caption{Explicit-solution validation.}
\label{tab:explicit_method_comparison}
\resizebox{\linewidth}{!}{%
\begin{tabular}{lrrrr}
\hline
Method
& $E_{x,\infty}$
& $E_{u,\infty}$
& $R_\infty$
& Time (s)\\
\hline
Discrete
& $2.42\times10^{-2}$
& $3.84$
& $3.80$
& 0.019\\
Bernstein
& $6.64\times10^{-4}$
& $7.19\times10^{-2}$
& $9.93\times10^{-7}$
& 0.979\\
Chebyshev
& $9.49\times10^{-4}$
& $1.60\times10^{-1}$
& $9.62\times10^{-7}$
& 0.960\\
Legendre
& $9.49\times10^{-4}$
& $1.60\times10^{-1}$
& $9.62\times10^{-7}$
& 1.080\\
Fourier--polynomial
& $4.69\times10^{-1}$
& $3.95$
& $1.02\times10^{-6}$
& 11.058\\
Mixed Bernstein--Fourier
& $2.85\times10^{-4}$
& $1.32\times10^{-2}$
& $9.92\times10^{-7}$
& 1.229\\
KL
& $2.08\times10^{-3}$
& $5.11\times10^{-3}$
& $1.16\times10^{-1}$
& 21.922\\
FE
& $1.19\times10^{-3}$
& $9.24\times10^{-2}$
& $9.59\times10^{-7}$
& 3.285\\
\hline
\end{tabular}}

\end{table}

Table~\ref{tab:explicit_method_comparison} shows that
FE reproduces the optimal state with maximum error
$1.19\times10^{-3}$ while maintaining a sampled
dynamics residual below $10^{-6}$.
The mixed method gives smaller state and control
errors in this experiment.
The comparison also illustrates that a small
dynamics residual does not by itself imply a small
error relative to the optimal solution.

\FloatBarrier
\subsection{High-Frequency Tracking with a Kinematic
Bicycle Model}
\label{subsec:bicycle_high_frequency}

This experiment examines whether the Fourier
representation can describe an oscillatory tracking
task with fewer optimization variables than the
tested polynomial representations.

Let
\[
x=(p_x,p_y,\theta,v)^\top,
\qquad
u=(a,\delta_s)^\top,
\]
where $(p_x,p_y)$ denotes position, $\theta$ heading,
$v$ speed, $a$ acceleration, and $\delta_s$ steering
angle. The kinematic bicycle dynamics are
\begin{equation}
\begin{aligned}
\dot p_x &= v\cos\theta,\\
\dot p_y &= v\sin\theta,\\
\dot\theta &= \frac{v}{L}\tan\delta_s,\\
\dot v &= a.
\end{aligned}
\label{eq:num_bicycle_dynamics}
\end{equation}

We minimize the tracking objective
\begin{equation}
\begin{aligned}
J(x,u)
={}&
\int_0^T
\left[
(x-x_r)^\top Q(x-x_r)+u^\top Ru
\right]\,dt\\
&+
(x(T)-x_d)^\top Q_f(x(T)-x_d),
\end{aligned}
\label{eq:num_bicycle_cost}
\end{equation}
subject to \eqref{eq:num_bicycle_dynamics},
\begin{equation}
x(0)=(0,0,0,2)^\top,
\end{equation}
and
\begin{equation}
0\leq v\leq3,\qquad
|a|\leq2,\qquad
|\delta_s|\leq0.8.
\end{equation}
Here $T=10$, $L=1$, and
$x_d=(20,0,0,2)^\top$.
The terminal state is penalized, rather than imposed
as an equality constraint.

With $\omega=2\pi/T$, define
\begin{equation}
h(t)
=
0.03[1-\cos(12\omega t)]
+
0.012[1-\cos(20\omega t)].
\end{equation}
The desired reference is
\begin{equation}
x_r(t)
=
[2t,
h(t),
\arctan\!\bigl(\dot h(t)/2\bigr),
\sqrt{4+\dot h(t)^2}\,\,]^T.
\label{eq:num_bicycle_reference}
\end{equation}
Its lateral component contains $12$ and $20$
oscillations over the horizon.
The reference specifies the tracking task; it is
not an exact optimal solution.

The weights are
\begin{align}
Q&=\operatorname{diag}(10,1000,50,10),\\
R&=\operatorname{diag}(0.1,0.05),\\
Q_f&=\operatorname{diag}(100,1000,100,50).
\end{align}
The parameterized formulations enforce constraints
at $3000$ uniform points, with dynamics tolerance
$\delta=10^{-3}$.
In this experiment, Chebyshev and Legendre use modal
representations evaluated on the same uniform grid.

Costs and residuals are checked independently on
a dense validation grid.

\begin{table}[htbp]
\centering
\caption{High-frequency bicycle tracking.}
\label{tab:bicycle_high_frequency}
\resizebox{\linewidth}{!}{%
\begin{tabular}{llrrrr}
\hline
Method & Parameters & Variables & $J$ & $R_\infty$
& Time (s)\\
\hline
Chebyshev
& $d=80$ & 486 & 2.167488
& $2.07\times10^{-3}$ & 133.566\\
Legendre
& $d=80$ & 486 & 2.167488
& $2.07\times10^{-3}$ & 114.990\\
FE
& $N=34,\ q=4$ & 414 & 2.180732
& $1.0012\times10^{-3}$ & 96.492\\
Fourier--polynomial$^\dagger$
& $N=16$ & --- & --- & --- & ---\\
Bernstein$^\dagger$
& $d=40$ & --- & --- & --- & ---\\
Mixed Bernstein--Fourier$^\dagger$
& $B=12,\ F=20$ & --- & --- & --- & ---\\
\hline
\end{tabular}}

\smallskip
\begin{minipage}{\linewidth}
\footnotesize
$B$ and $F$ denote the Bernstein degree and Fourier
order of the mixed method.
$^\dagger$The recorded trial ended with an
implementation exception (\texttt{SystemError});
no usable numerical result or solver time was
recorded. These entries are not classified as
time-limit terminations.
\end{minipage}
\end{table}

FE obtains a cost approximately $0.61\%$ above the
Chebyshev and Legendre values while using $14.8\%$
fewer optimization variables.
It also has a shorter recorded solution time and
approximately half their sampled dynamics residual.
Its residual slightly exceeds the nominal
$10^{-3}$ tolerance on the independent grid.

These results support the efficiency of the tested
FE representation for this oscillatory problem.
They do not establish a general superiority over
polynomial methods at all approximation orders.

\begin{figure}[htbp]
\centering
\includegraphics[width=\linewidth]{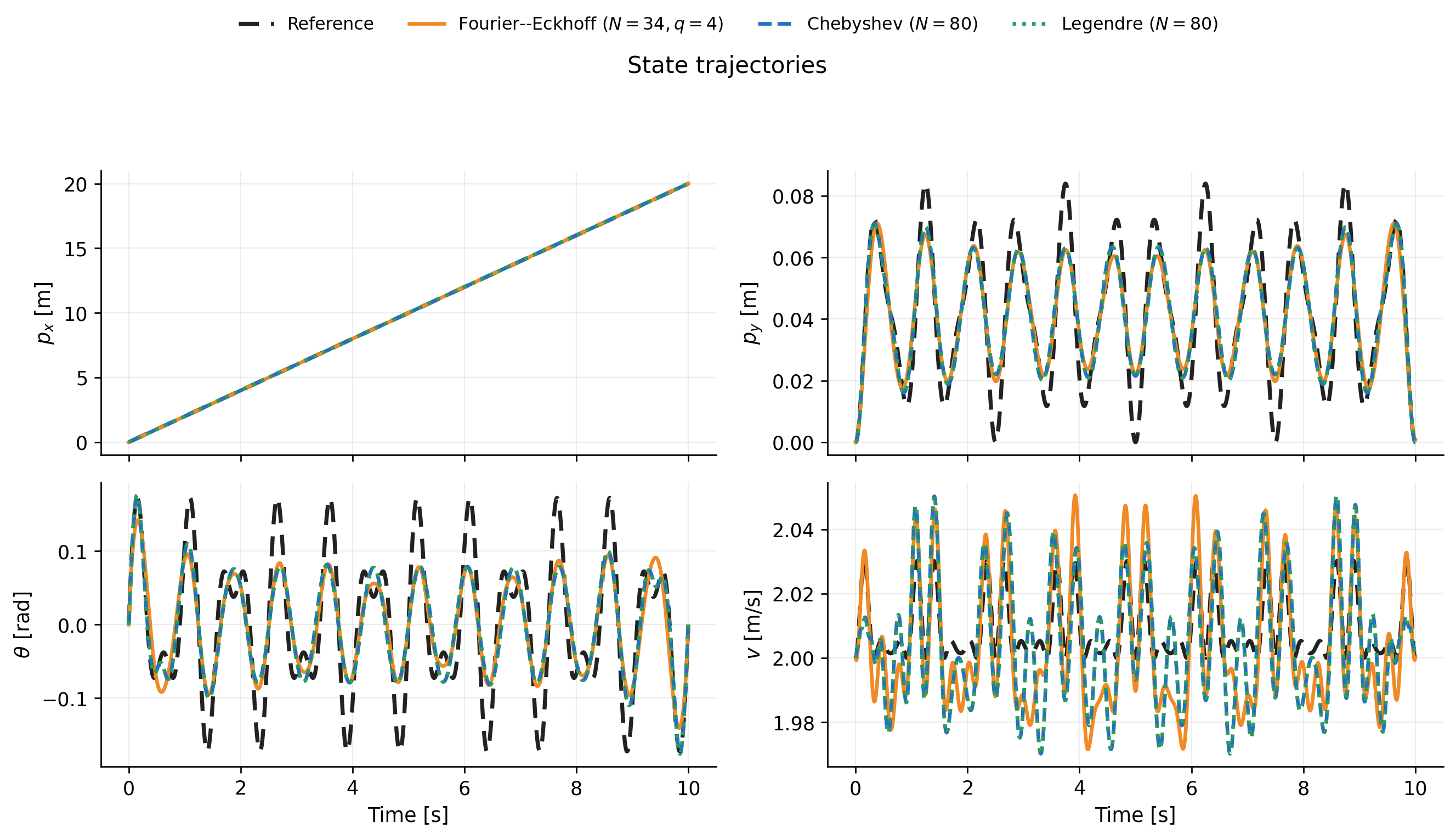}
\caption{Bicycle states.}
\label{fig:bicycle_states}
\end{figure}

\begin{figure}[htbp]
\centering
\includegraphics[width=\linewidth]{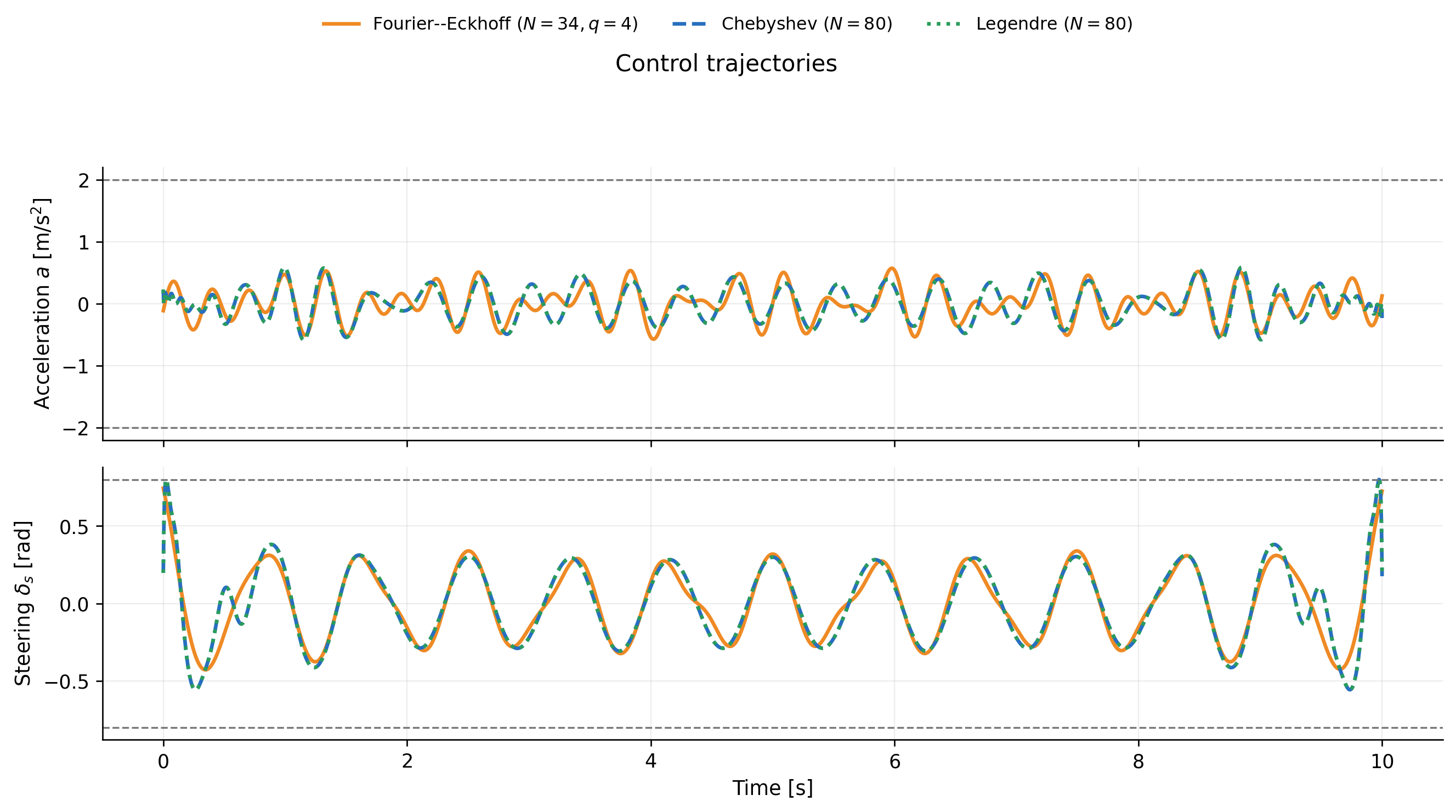}
\caption{Bicycle controls.}
\label{fig:bicycle_controls}
\end{figure}

\FloatBarrier
\subsection{Uniform Grids for Fixed-Rate Sensing}
\label{subsec:fixed_rate_sensing}

The purpose of this experiment is to examine an
application in which uniform sampling is imposed
by the measurement process.
A moving observer selects its trajectory to improve
bearing-only information about a stationary target,
while satisfying motion and obstacle constraints.

The observer state and control are
\[
z=(p_x,p_y,\psi)^\top,
\qquad
w=(v,\omega_u)^\top,
\]
with dynamics
\begin{equation}
\dot p_x=v\cos\psi,\qquad
\dot p_y=v\sin\psi,\qquad
\dot\psi=\omega_u.
\label{eq:num_observer_dynamics}
\end{equation}
We use $T=20$ seconds and impose
\begin{equation}
z(0)=(-10,-2,0)^\top,\qquad
z(T)=(10,2,0)^\top.
\end{equation}
The path constraints are
\begin{equation}
0.2\leq v\leq2.2,\qquad
|\omega_u|\leq1.8,\qquad
|p_y|\leq3.
\end{equation}

The target is at $\xi=(0,0)^\top$.
The observer must satisfy
$\|p(t)-\xi\|_2\geq1.10$, where
$p=(p_x,p_y)^\top$.
Six circular obstacles have centers
\begin{equation}
c_j=
\begin{bmatrix}
-7+3(j-1)\\
0.65(-1)^{j+1}
\end{bmatrix},
\qquad j=1,\ldots,6.
\end{equation}
Their physical radius is $0.58$, and the required
clearance is $0.18$. Thus,
\begin{equation*}
\|p(t)-c_j\|_2\geq0.76,
\qquad j=1,\ldots,6.
\end{equation*}

Measurements arrive at the fixed rate $f_s=25$ Hz:
\begin{equation*}
t_k=\frac{k}{f_s},\qquad k=0,\ldots,500.
\end{equation*}
Let $d_k=\xi-p(t_k)$ and define
\begin{equation}
g_k=
\frac{1}{\|d_k\|_2^2+10^{-8}}
\begin{bmatrix}
-d_{k,y}\\ d_{k,x}
\end{bmatrix}.
\end{equation}
The normalized information matrix is
\begin{equation*}
\mathcal I
=
0.05I_2+
\frac{1}{501\,\sigma^2}
\sum_{k=0}^{500}g_kg_k^\top,
\qquad \sigma=0.08.
\label{eq:num_information_matrix}
\end{equation*}

The optimal control problem minimizes
\begin{equation*}
J_{\mathrm{obs}}
=
-\log\det\mathcal I
+
0.025\int_0^{20}
\left[(v-1.35)^2+0.25\omega_u^2\right]\,dt
\label{eq:num_observer_cost}
\end{equation*}
subject to the dynamics, endpoint conditions,
and path constraints above.
A larger $\log\det\mathcal I$ indicates greater
information, although the optimized objective also
includes the motion penalty.

All methods use the same $501$ sensor times for
information accumulation and sampled path
constraints.
FE and Bernstein also enforce dynamics at these
times. The tested Chebyshev and Legendre
formulations enforce dynamics at their respective
$d+1$ Lobatto nodes.
The dynamics tolerance in the optimization is
$8\times10^{-4}$, and validation uses $25\,001$
uniform points.

\begin{table}[htbp]
\centering
\caption{Fixed-rate sensing with six obstacles.}
\label{tab:fixed_rate_sensing}
\resizebox{\linewidth}{!}{%
\begin{tabular}{llrrrr}
\hline
Method & Parameters & Variables
& $\log\det\mathcal I$ & $R_\infty$ & Time (s)\\
\hline
Bernstein
& $d=24$ & 125 & 7.18622
& $9.13\times10^{-4}$ & 114.12\\
FE
& $N=12,\ q=4$ & 125 & 7.21268
& $8.06\times10^{-4}$ & 91.55\\
Chebyshev
& $d=40$ & 205 & 7.32487
& $2.08\times10^{-1}$ & 141.58\\
Legendre
& $d=40$ & 205 & 7.32048
& $1.82\times10^{-1}$ & 131.12\\
FE
& $N=20,\ q=4$ & 205 & 7.30996
& $8.21\times10^{-4}$ & 271.44\\
\hline
\end{tabular}}
\end{table}

At $125$ optimization variables, FE provides slightly
greater information than Bernstein, with a smaller
dynamics residual and shorter recorded computing
time.
At $205$ variables, FE provides information close
to the Chebyshev and Legendre results while reducing
the sampled dynamics residual by more than two
orders of magnitude. This improvement requires
more dynamics-enforcement points and a longer
solution time.

The experiment illustrates the convenience of
aligning the FE discretization with physical sensor
times. However, it does not isolate grid uniformity
from the number of enforcement points.
Polynomial representations can also be evaluated
on uniform grids, as in the preceding experiment.
Furthermore, small negative obstacle-clearance
values were observed during dense validation;
sampled enforcement alone therefore does not
establish continuous-time obstacle feasibility.

\begin{figure}[htbp]
\centering
\includegraphics[width=\linewidth]{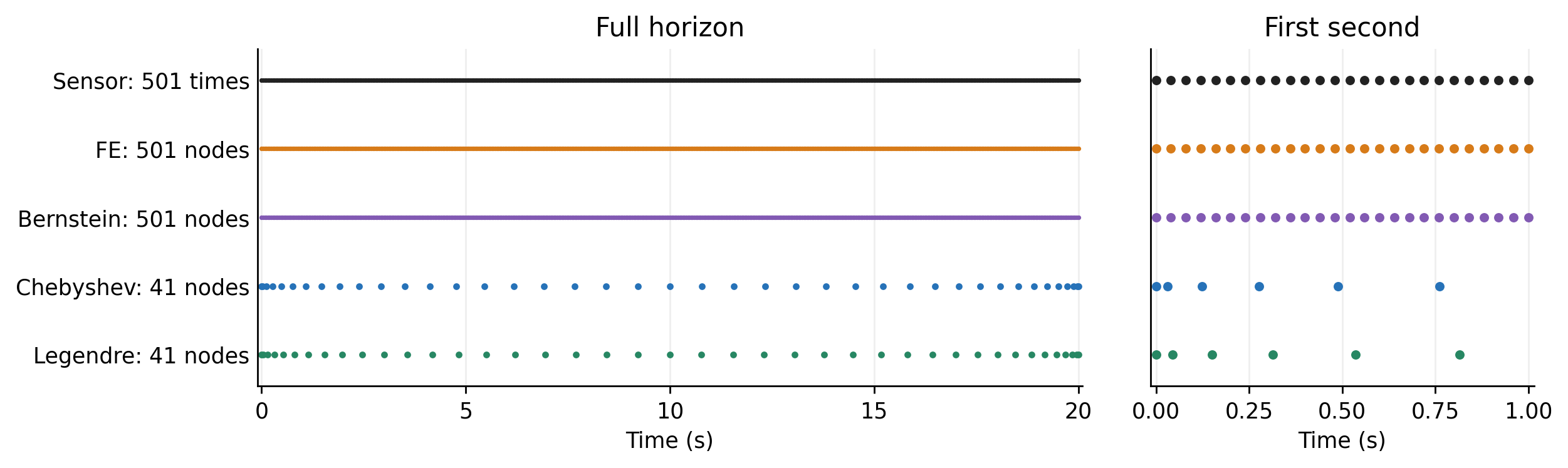}
\caption{Sensor times and dynamics-enforcement grids.}
\label{fig:sensor_times}
\end{figure}

\FloatBarrier
\section{Sobolev Constraints and Five-Obstacle
Avoidance}
\label{subsec:sobolev_obstacle_avoidance}

The final experiment addresses a limitation of
pointwise constraint enforcement: a reconstructed
trajectory may satisfy a constraint at the
optimization nodes and violate it between them.
We first describe the Sobolev formulation and then
evaluate it on a five-obstacle avoidance problem.

\subsubsection{Controlling the dynamics residual}

Assume that $r_N\in H^1(0,T;\mathbb R^{n_x})$ and define
\begin{equation*}
\|r_N\|_{H^1}^2
=
\sum_{i=1}^{n_x}
\int_0^T
\left[r_{N,i}(t)^2+\dot r_{N,i}(t)^2\right]\,dt.
\end{equation*}
The one-dimensional Sobolev estimate gives
\begin{equation*}
\sup_{t\in[0,T]}\|r_N(t)\|_\infty
\leq
C_T\|r_N\|_{H^1},
\qquad
C_T=\sqrt{1+\frac{1}{T}}.
\label{eq:num_sobolev_embedding}
\end{equation*}

Hence the single scalar constraint
\begin{equation*}
\|r_N\|_{H^1}^2
\leq
\left(\frac{\delta}{C_T}\right)^2
\label{eq:num_sobolev_dynamics}
\end{equation*}
is sufficient to bound every dynamics-residual
component by $\delta$ throughout the time interval.

In the numerical implementation, let
$t_k=k\Delta t$, where
$\Delta t=T/(K-1)$.
For a scalar function $v$, define
\begin{equation}
\begin{aligned}
\mathcal H_K(v)^2
={}&
\sum_{k=0}^{K-1}w_kv(t_k)^2\\
&+
\frac{1}{\Delta t}
\sum_{k=0}^{K-2}
\left[v(t_{k+1})-v(t_k)\right]^2,
\end{aligned}
\label{eq:num_discrete_h1}
\end{equation}
where
\[
w_0=w_{K-1}=\frac{\Delta t}{2},
\qquad
w_k=\Delta t,\quad 1\leq k\leq K-2.
\]
For a vector residual, the squared quantities are
summed over its components.
The implemented dynamics constraint is
\begin{equation}
\mathcal H_K(r_N)^2
\leq
\left(\frac{\delta}{C_T}\right)^2.
\label{eq:num_discrete_sobolev_dynamics}
\end{equation}

It is essential to distinguish
\eqref{eq:num_discrete_sobolev_dynamics}
from the exact norm constraint
\eqref{eq:num_sobolev_dynamics}.
If a verified error bound satisfies
\begin{equation*}
\|r_N\|_{H^1}
\leq
\mathcal H_K(r_N)+\eta_K,
\end{equation*}
then the tightened condition
\begin{equation*}
\mathcal H_K(r_N)+\eta_K
\leq\frac{\delta}{C_T}
\end{equation*}
guarantees the continuous-time residual bound.
Without such an error bound, dense-grid checks
provide numerical evidence rather than a rigorous
continuous-time certificate.

\subsubsection{Extending the idea to obstacle constraints}

A bound on the dynamics residual does not directly
control obstacle clearance.
To treat an obstacle inequality
\begin{equation*}
h_j(x_N(t))\geq0,
\end{equation*}

we introduce a scalar slack trajectory $s_j(t)$
and the residual
\begin{equation}
e_j(t)
=
h_j(x_N(t))-s_j(t)^2-\rho,
\qquad \rho>0.
\label{eq:num_safety_residual}
\end{equation}
If $e_j\in H^1(0,T)$ and
\begin{equation*}
\|e_j\|_{H^1}\leq\varepsilon_s,
\qquad
\rho\geq C_T\varepsilon_s,
\label{eq:num_safety_sobolev}
\end{equation*}
then, for every $t\in[0,T]$,
\begin{equation}
\begin{aligned}
h_j(x_N(t))
&=s_j(t)^2+\rho+e_j(t)\\
&\geq \rho-|e_j(t)|\\
&\geq \rho-C_T\varepsilon_s
\geq0.
\end{aligned}
\label{eq:num_safety_bound}
\end{equation}
Thus, an exact Sobolev bound on each safety residual
implies obstacle avoidance throughout the interval.

The slack trajectories are represented by FE
expansions. In the implementation, their safety
residual norms are approximated using
\eqref{eq:num_discrete_h1}.
The same distinction between an exact norm bound
and its numerical approximation applies here.

\subsubsection{Five-obstacle optimal control problem}

The purpose of this experiment is to compare
pointwise constraint enforcement with simultaneous
Sobolev control of dynamics and obstacle residuals.

We use the bicycle dynamics
\eqref{eq:num_bicycle_dynamics} and objective
\eqref{eq:num_bicycle_cost}, now with
\begin{equation*}
T=10,\qquad L=1.5,
\end{equation*}
\begin{equation*}
x(0)=(0,0,0,2)^\top,\qquad
x_d=(20,0,0,2)^\top,
\end{equation*}
and the straight-line reference
\begin{equation*}
x_r(t)=(2t,0,0,2)^\top.
\end{equation*}
The terminal state remains a soft target.
The weights are
\begin{align*}
Q&=\operatorname{diag}(5,5,2,2),\\
R&=\operatorname{diag}(0.2,0.1),\\
Q_f&=\operatorname{diag}(200,300,80,40).
\end{align*}
The bounds are
\begin{equation*}
0.5\leq v\leq3,\qquad
|a|\leq2,\qquad
|\delta_s|\leq0.7,\qquad
|p_y|\leq3.
\end{equation*}

Five obstacles are placed at
\begin{equation*}
c_j=
\begin{bmatrix}
4+3(j-1)\\
0.48(-1)^{j+1}
\end{bmatrix},
\qquad j=1,\ldots,5.
\label{eq:num_five_obstacle_centers}
\end{equation*}
Each obstacle has physical radius $0.54$, prescribed
clearance $0.12$, and additional margin $0.05$.
The inflated safety radius is therefore
\begin{equation*}
R_{\mathrm{safe}}=0.54+0.12+0.05=0.71.
\end{equation*}
The obstacle constraints are
\begin{equation*}
h_j(x(t))
=
\|p(t)-c_j\|_2^2-R_{\mathrm{safe}}^2
\geq0,
\qquad j=1,\ldots,5.
\label{eq:num_five_obstacle_constraints}
\end{equation*}

The straight-line reference intersects the inflated
obstacles, so a feasible trajectory must perform
repeated lateral avoidance maneuvers.

The parameterized formulations use $K=601$ uniform
nodes and dynamics tolerance $\delta=10^{-3}$.
The discrete baseline uses $K=1001$ nodes.
All methods are initialized from the same slalom
construction, represented in their respective
approximation spaces.

The pointwise formulations impose the dynamics and
obstacle inequalities at the optimization nodes.
FE-S instead uses one aggregate Sobolev constraint
for the dynamics and one Sobolev constraint for each
of the five safety residuals.
The remaining state and control bounds are enforced
at the nodes.

For FE-S, the state and control order is $N=24$,
the slack order is $N_s=72$, and $q=4$.
We use
\begin{equation*}
\rho=0.006,\qquad
\varepsilon_s=0.005.
\end{equation*}
Since $C_T=\sqrt{1.1}$,
\[
\rho>C_T\varepsilon_s.
\]
The implemented safety-norm constraints include
additional tightening:
\begin{equation*}
\mathcal H_K(e_j)
\leq0.9\varepsilon_s=0.0045,
\qquad j=1,\ldots,5.
\end{equation*}
This tightening provides numerical margin but is
not, by itself, a verified bound on norm
approximation error.

Validation uses $100\,001$ uniform points.
Basis representations are differentiated
analytically; the discrete solution is reconstructed
piecewise linearly.
In addition to cost and dynamics residual, we report
the minimum sampled safety margin
\begin{equation}
m_{\min}
=
\min_{\substack{t\in\mathcal G_{\mathrm{val}}\\
                j=1,\ldots,5}}
\left(\|p_N(t)-c_j\|_2-R_{\mathrm{safe}}\right).
\label{eq:num_minimum_clearance}
\end{equation}
A negative value indicates a violation of the
inflated safety boundary, not necessarily an
intersection with the physical obstacle.

\begin{table}[htbp]
\centering
\caption{Five-obstacle avoidance.}
\label{tab:five_obstacle_comparison}
\resizebox{\linewidth}{!}{%
\begin{tabular}{llrrrrr}
\hline
Method & Parameters & Variables & $J$
& $R_\infty$ & $m_{\min}$ (mm) & Time (s)\\
\hline
Discrete
& $K=1001$
& 6006 & 1.205274
& $4.31\times10^{-1}$ & $-0.0311$ & 0.687\\
Bernstein
& $d=24$
& 150 & 1.785481
& $1.02\times10^{-3}$ & $-0.0228$ & 29.473\\
Chebyshev
& $d=40$
& 246 & 1.432037
& $1.09\times10^{-3}$ & $-0.0448$ & 27.080\\
Legendre
& $d=40$
& 246 & 1.432037
& $1.09\times10^{-3}$ & $-0.0448$ & 26.632\\
Fourier--polynomial
& $N=18$
& 240 & 1.344416
& $1.02\times10^{-3}$ & $-0.0242$ & 195.261\\
FE
& $N=24,\ q=4$
& 294 & 1.292812
& $1.02\times10^{-3}$ & $-0.0300$ & 28.424\\
FE-S$^\dagger$
& $N=24,\ N_s=72,\ q=4$
& 1019 & 1.731065
& $3.33\times10^{-5}$ & $+1.9776$ & 209.815\\
Mixed Bernstein--Fourier$^\ddagger$
& $B=6,\ F=12$
& 192 & 1.454210
& $1.03\times10^{-3}$ & $-0.0790$ & 281.366\\
\hline
\end{tabular}}

\smallskip
\begin{minipage}{\linewidth}
\footnotesize
$^\dagger$FE-S terminated at IPOPT's acceptable level.
Its variable count includes five slack trajectories.

$^\ddagger$The mixed-method row reports an
unsuccessful final iterate after
\texttt{Maximum\_CpuTime\_Exceeded}.
All four tested mixed configurations reached the
CPU-time limit.
The configured IPOPT CPU-time cap was $270$ s for
Bernstein, Fourier--polynomial, and mixed methods,
and $240$ s for the other methods.
Reported solver times are elapsed wall times and
can exceed these CPU-time caps.
\end{minipage}
\end{table}

Among the converged pointwise parameterized
formulations in Table~\ref{tab:five_obstacle_comparison},
FE gives the lowest recorded cost, with computing
time close to those of Bernstein, Chebyshev,
and Legendre.
However, all these pointwise solutions exhibit
small negative safety margins on the dense
validation grid.
The discrete baseline has a still lower cost,
but its large reconstructed dynamics residual
prevents interpreting this as a comparable
dynamically feasible result.

FE-S reduces the sampled dynamics residual by
approximately a factor of $30$ relative to FE and
maintains a positive minimum sampled safety margin
of approximately $1.98$ mm.

The independently evaluated maximum safety-residual
$H^1$ norm is approximately $4.60\times10^{-3}$,
below the nominal value $\varepsilon_s=0.005$.

These improvements require additional slack
variables, a higher objective value, and a longer
solution time.

\begin{figure}[htbp]
\centering
\includegraphics[width=\linewidth]{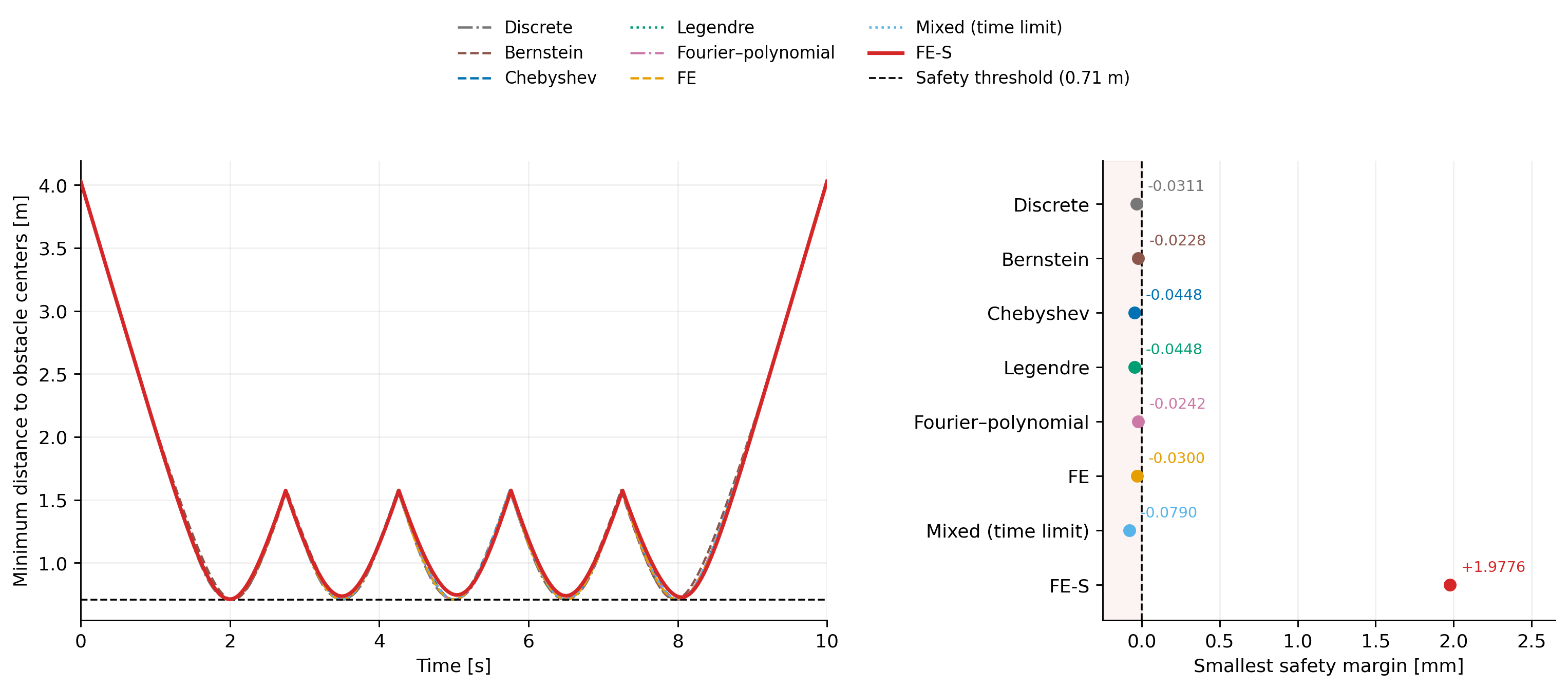}
\caption{Obstacle distance and minimum safety margin for Table~\ref{tab:five_obstacle_comparison}.}
\label{fig:five_obstacle_distance}
\end{figure}

Figure~\ref{fig:five_obstacle_distance} compares the minimum obstacle
distances for the configurations in Table~\ref{tab:five_obstacle_comparison}.
Increasing the approximation order alone did not ensure satisfaction
of the prescribed safety margin in the additional recorded runs.
For Chebyshev and Legendre, increasing the degree from $40$ to $56$
reduced the minimum-margin violation from approximately $0.0448$ mm
to $0.0034$ mm, but did not eliminate it. Bernstein runs at degrees
$28$ and $32$ also retained negative safety margins and reached the
$270$ s CPU-time limit. These are violations of the inflated safety
boundary, not observed collisions with the physical obstacles.
Higher order increases the number of decision variables and can
increase computing time, although the measured times are not
monotone in degree. Thus, FE-S's longer solution time buys explicit
Sobolev residual control and a positive sampled safety margin in this
experiment; increasing the order of a pointwise formulation alone
does not provide the same guarantee. The continuous-time safety
guarantee in \eqref{eq:num_safety_bound} remains conditional on
satisfaction of the stated Sobolev norm bound.

The experiment therefore demonstrates a trade-off:
the pointwise FE formulation provides a compact
and comparatively inexpensive approximation,
whereas the Sobolev extension gives stronger
numerical control of dynamics and obstacle
residuals.
The continuous-time implication
\eqref{eq:num_safety_bound} is mathematical;
the reported discretized computations support it
numerically but do not supply a verified
quadrature-error certificate.

\begin{figure}[htbp]
\centering
\includegraphics[width=\linewidth]{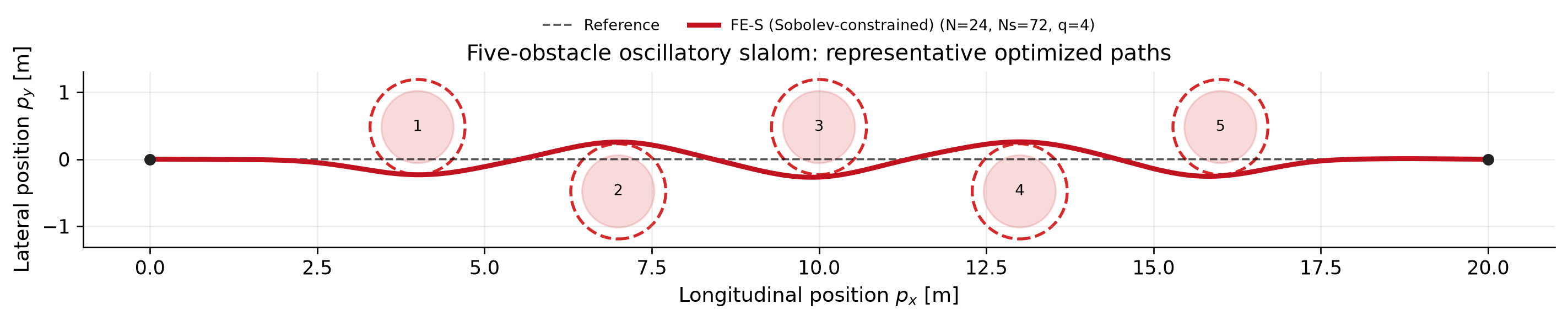}
\caption{FE-S path around five obstacles.}
\label{fig:five_obstacle_path}
\end{figure}

In the five-obstacle experiment, FE-S reduces the sampled dynamics residual by approximately a factor of $30$ and yields a positive sampled safety margin. These improvements require additional variables, a higher objective value, and more computing time. The discretized results support improved residual control, while continuous-time guarantees remain conditional on verified norm bounds.

\bibliographystyle{plain}
\bibliography{references}

\end{document}